
\input amstex
\documentstyle{amsppt}
\magnification 1200
\NoBlackBoxes
\nologo

\define\a{\alpha}
\define\be{\beta}
\define\ka{\kappa}

\define\la{\lambda}
\define\U{U_q(\widehat{\goth g})}
\define\alg{U_q(\widehat{\goth g})}
\define\gA{\goth g}
\define\Z{Z(\epsilon_1\cdots\epsilon_s;z_1\cdots z_s)}
\define\W{Z(\epsilon_1\cdots\epsilon_s;n_1\cdots n_s)}
\define\e{\epsilon}
\redefine\v{Z(\epsilon_1\cdots\epsilon_s;z_1\cdots z_s)v_{\lambda}}
\redefine\u{Z(\epsilon_1\cdots\epsilon_s;n_1\cdots n_s)v_{\lambda}}
\define\vep{\varepsilon}

\document

\topmatter
\title Quantum Z-algebras and representations of quantum affine algebras
\endtitle
\rightheadtext{Quantum $Z$-algebras}
\author Naihuan Jing\endauthor
\affil Department of Mathematics\\
       North Carolina State University\\
       Raleigh, NC 27695-8205, U.S.A.
\endaffil
\email jing\@eos.ncsu.edu\endemail
\subjclass 17B\endsubjclass
\abstract 
Generalizing our earlier work, we introduce the homogeneous
quantum $Z$-algebras for all quantum affine algebras $\alg$ of type one. 
With the new algebras we unite previously scattered realizations
of quantum affine algebras in various cases. As a result we find a 
realization of $U_q(F_4^{(1)})$.
\endabstract

\thanks Supported in part by NSA grants MDA904-96-1-0087 and MDA904-96-1-0087.
\endthanks
\keywords Quantum enveloping algebras,
Quantum $Z$-operators, Vertex operators, Hopf algebras\endkeywords
\endtopmatter

\head 0. Introduction
\endhead

In 1981 Lepowsky and Wilson introduced (principal)
$Z$-algebras as a tool to construct
explicit principal representations of affine Lie algebras \cite{LW1}.
Later Lepowsky and Primc considered the homogeneous analog for $\widehat{sl}_2$
and gave explicit realizations of higher level representations of
$\widehat{sl}_2$ \cite{LP1, LP2}. 
Such a construction is also done for the higher rank cases
\cite{P}. Many other vertex representations of affine algebras were constructed
by this method (cf. \cite{FLM} for
a list of references).

Recently we introduced the quantum $Z$-algebras for the simply laced
types in the homogeneous picture \cite{J2} and used them to construct
higher level (integral) representations of the quantum affine algebra
$U_q(\widehat{sl}_2)$. In \cite{BV} the higher level representations of
$U_q(\widehat{sl}_2)$ and its $Z$-algebra
in the Verma module picture were also studied.
These constructions parallel to the
classical method of Lepowsky and Primc. As we mentioned in \cite{J2},
the quantum $Z$-algebras are expected to provide more examples of explicit
constructions of quantum affine algebras. 

In this paper we define the quantum $Z$-algebras for all type
one  quantum affine algebras. The quantum $Z$-algebras enable us to unify 
previously 
scattered constructions of quantum affine algebras, and we also find 
a realization of $U_q(F_4^{(1)})$ based on our recent joint work on 
representations of $U_q(C_n^{(1)})$ \cite{JKM1, JKM2}. The latter also 
contained some
quantum deformation of the admissible representations 
studied by Kac and Wakimoto in the classical
case (in this special case, see Feingold-Frenkel \cite{FF}). Thus all level one
representations of quantum affine algebras of type one can be constructed by
vertex operators and quantum $Z$-algebras.
Since we work in a general setting
 we can derive all commutation relations for the 
quantum $Z$-algebras at once with the help of $q$-series.
As an simple example we recover the Frenkel-Jing bosonic realizations
of quantum affine algebras of simply laced types. As for other cases
we try to emphasize on the 
similarity among all the constructions and points out differences in 
each special case.

This paper mainly concerns with homogeneous 
quantum $Z$-algebras. The quantum 
principal $Z$-algebra is still unknown, though the related
fermionic constructions are known quite generally (cf. \cite{KMPY}), where
one can also find more examples about higher level 
$U_q(\widehat{sl_2})$-modules.

\head 1. Quantum Affine Algebras and their Modules
\endhead

Let $\widehat{\goth g}$ be the affine Kac-Moody algebra \cite{K} 
associated to generalized Cartan matrix of 
type one $A=(A_{ij})$. Let $\gA$ be the associated finite dimensional 
simple Lie algebra. Let $Q=\Bbb Z\a_1\oplus \cdots \oplus
\Bbb Z\a_l$ and $P=\Bbb Z\la_1\oplus \cdots \oplus\Bbb Z\la_l$ be the 
root and weight lattices of $\gA$. We fix the standard normalized invariant
$(\ |\ )$ on $\gA$ such that
$$(\a_i|\a_j)=d_iA_{ij}, \qquad i, j=1, \cdots l
$$
where the $d_i=(\a_i|\a_i)/2$. Let $q_i=q^{d_i}$.

The quantum affine algebra $\alg$  is the associative algebra 
generated by 
the central element $\gamma$, the grading element $d$, and the 
elements $\a_i(n)$, and $x_i^{\pm}(n)$, $i=1, \cdots l$ (or equivalently the elements
$x^{\pm}(n), \phi(-m), \psi(m)$), $n\in \Bbb Z$ and $m\in \Bbb Z_+$ subject to
the following relations:
$$\gather
[d, a_i(n)]=n\a_i(n), [d, x^{\pm}_i(n)]=nx_i^{\pm}(n),\\
[a_i(m), a_j(n)]=\delta_{m,-n}\frac{[(\a_i|\a_j)m]}{m}
\frac{\gamma^m-\gamma^{-m}}{q-q^{-1}},\\
[a_i(m), x_j^{\pm}(n)] =\pm \frac{[(\a_i|\a_j)m]}{m}\gamma^{\mp |m|/2}
x_j^{\pm}(m+n),\\
(z-q^{\pm A_{ij}}w)x_i^{\pm}(z)x_j^{\pm}(w)=x_j^{\pm}(w)x_i^{\pm}(z)(q^{\pm
A_{ij}}z-w)\\
[x^+_i(z), x^-_j(w)]=\frac{\delta_{ij}}{q_i-q_i^{-1}}
\{\delta(\frac zwq^{-c})\psi_i(wq^{c/2})-\delta(\frac zwq^c)
\phi_i(zq^{c/2})\},\\
Sym_{z_1, \cdots z_N}\sum_{s=0}^{N=1-A_{ij}}(-1)^s\bmatrix N\\s\endbmatrix_i
x^{\pm}_i(z_1)\cdots x_i^{\pm}(z_s)\cdot\\
{\kern .5cm}\cdot x^{\pm}_j(w)x_i^{\pm}(z_{s+1})\cdots
x^{\pm}_i(z_N)=0, \quad\text{for}\quad A_{ij}\neq 0
\endgather
$$
where the generators $\a(n)$ are related to $\phi(-m), \psi(m)$ via:
$$\gather
\phi_i(z)=\sum_{n\geq 0}\phi_i(-n)z^n=k_i^{-1}exp((q^{-1}-q)\sum_{n>0}\a_i(-n)z^n),\\
\psi_i(z)=\sum_{n\geq 0}\psi_i(n)z^{-n}=k_iexp((q-q^{-1})\sum_{n>0}\a_i(n)z^{-n}).
\endgather
$$
The fields are defined by
$$
x_i^{\pm}(z)=\sum_{n\in \Bbb Z}x_i^{\pm}(n)z^{-n-1}, \qquad i=1, \cdots, l.
$$
We have rescaled the generators $a_i(n)$ by the factor $[(\a_i|\a_i)/2]$ to 
have an invariant presentation. With $a_i(n)/[(\a_i|\a_i)/2]$ as $a_i(n)$ 
we get the Drinfeld realization studied in \cite{Be, J4}.
 
We define the $q$-integers by
$$
[n]_i=\frac{q_i^n-q_i^{-n}}{q_i-q_i^{-1}}, \quad
\bmatrix m\\n\endbmatrix _i=\frac{[m]_i!}{[n]_i![m-n]_i!}.
$$
We also need the $q$-number (for $a\in \Bbb R$):
$$(a;q)_n=(1-a)(1-aq)\cdots (1-q^{n-1}a).
$$

Let $V$ be a $\alg$-module. The module $V$ is called level $k$
if
the central element $\gamma$ acts as the scalar $q^k$ on $V$.
For $a\in \Bbb C$, define
$$V_a=\{v\in V| d.v=av\},
$$
the elements of $V_a$ are homogeneous elements of degree $a$. An endomorphism 
$A$ of
$V$ is homogeneous of degree $n\in \Bbb C$ if 
$$
[d, A]=nA.
$$

We will study the following category of $\U$-modules $V$ with the properties:
\roster
\item the module $V$ has level $k$;
\item each weight space of $V$ is finite dimensional;
\item for each $a\in \Bbb C$, there exists $n_0\geq 0$ such that for
all $n>n_0, V_{a+n}=\{0\}$.
\endroster

It is clear that the irreducible highest weight modules $L(\lambda)$ 
lie in the category for a suitable $k$. The level is an integer for 
the dominant 
integrable module.

\head 2. Quantum $Z$-algebras
\endhead

Let $V$ be a $\alg$-module of level $k$, $k\in \Bbb Q$. 
The Heisenberg algebra $U_q(\widehat{\goth h})$ of level $k$
is generated by $\a_i(n), n\neq 0$ and the central element $c$ 
subject to the relations:
$$
[\a_i(m), \a_j(n)]=\frac 1m [(\a_i|\a_j)m][mk]\delta_{m,-n}.\tag2.1
$$
and $V$ is also a $U_q(\widehat{\goth h})$-module.

We define the following operators on $V$:
$$\aligned
E_-^{\pm}(\a_i, z)&=exp(\mp\sum_{n\geq 1}\frac{q^{\mp kn/2}}{[kn]}\a_i(-n)z^n)\\
E_+^{\pm}(\a_i, z)&=exp(\pm\sum_{n\geq 1}\frac{q^{\mp kn/2}}{[kn]}\a_i(n)z^{-n})
\endaligned\tag2.2
$$

For $a\in \Bbb R$ we set
$$
(1-z)_{q^{2k}}^{a/k}=\frac{(q^{-a+k}z;q^{2k})_{\infty}}
{(q^{a+k}z;q^{2k})_{\infty}}=exp(-\sum_{n=1}^{\infty}\frac{[an]}{n[kn]}z^n)
$$

The Heisenberg algebra has a natural realization on the space 
$Sym(\widehat{\goth h}^-)$ of the symmetric algebra generated by 
the elements $\a_i(-m)$, $m\in \Bbb N$ via the following rule:
$$\align
q^c.1&=q^k, \qquad\qquad q^{-c}.1=q^{-k},\\
\a_i(-n)&=\text{multiplication operator}, \quad n\in \Bbb N\\
\a_i(n)&=\text{differentiation operator subject to (2.1)}.
\endalign
$$
We denote the resulted irreducible module of $U_q(\widehat{\goth h})$ as $K(k)$.

We now introduce the $Z$-vertex operators. Let $V$ be a level $k$ 
$\U$-module. 
For each $i=1, \cdots, l$, define
$$
Z^{\pm}(\a_i, z)=E^{\pm}_-(\a_i,z)x^{\pm}_i(z)E^{\pm}_+(\a_i, z), \tag2.5
$$
which act on the space $V$.

\proclaim{Proposition 2.1} As operators on a level $k$ module $V$,
$$\align
E_+^{\pm}(\a_i,z)E_-^{\pm}(\a_j,w)&=E_-^{\pm}(\a_j, w)E_+^{\pm}(\a_i,z)
(1-\frac wz q^{\mp k})_{q^{2k}}^{(\a_i|\a_j)/k}, \\
E_+^{\pm}(\a_i,z)E_-^{\mp}(\a_j,w)&=E_-^{\mp}(\a_j, w)E_+^{\pm}(\a_i,z)
(1-\frac wz)_{q^{2k}}^{-(\a_i|\a_j)/k},\\
E_+^{\pm}(\a_i,z)x^{\pm}_j(w)&=x^{\pm}_j(w)E_+^{\pm}(\a_i,z)
(1-\frac wz q^{\mp k})_{q^{2k}}^{-(\a_i|\a_j)/k},\\
E_+^{\pm}(\a_i,z)x^{\mp}_j(w)&=x^{\mp}_j(w)E_-^{\pm}(\a_i,z)
(1-\frac wz)_{q^{2k}}^{(\a_i|\a_j)/k},\\
E_+^{\pm}(\a_i,z)x^{\pm}_j(w)&=x^{\pm}_j(w)E_+^{\pm}(\a_i,z)
(1-\frac wzq^{\mp k})^{-(\a_i|\a_j)/k},\\
[\phi(\a_i, z), E_-^{\pm}(\a_j,w)]&=[\psi(\a_i,z), E_+^{\pm}(\a_j,w)]=0,\\
E_+^{\pm}(\a_j,w)\phi(\a_i,z)&=\phi(\a_i,z)E_+^{\pm}(\a_j,w)
\frac{w-q^{\pm (\a_i|\a_j)\mp k/2}z}{w-q^{\mp (\a_i|\a_j)\mp k/2}z},\\
\psi(\a_i,z)E_-^{\pm}(\a_j,w)&=E_-^{\pm}(\a_j,w)\psi(\a_i,z)
\frac{z-q^{\pm (\a_i|\a_j)\mp k/2}w}{z-q^{\mp (\a_i|\a_j)\mp k/2}w}.\\
E^+_-(\a_i, z)E^-_-(\a_i, zq^{-k})&=1,\\
E^+_-(-\a_i, z)E^-_-(-\a_i, zq^k)&=\phi(\a_i, zq^{k/2})k_i,\\
E^+_+(-\a_i, z)E^-_+(-\a_i, zq^{-k})&=\psi(\a_i, zq^{-k/2})k_i^{-1},\\
E^+_+(\a_i, z)E^-_+(\a_i, zq^k)&=1.
\endalign
$$
\endproclaim \hfill$\square$

\proclaim{Theorem 2.2} As operators on a level $k$ space, 
the $Z$-vertex operators satisfy the following relations.
$$\align
&(1-q^{\mp k}w/z)_{q^{2k}}^{-(\a_i|\a_j)/k}(z-q^{\pm (\a_i|\a_j)}w)Z^{\pm}(\a_i,z)Z^{\pm}(\a_j,w)\\
&=(1-q^{\mp k}z/w)_{q^{2k}}^{-(\a_i|\a_j)/k}
(q^{\pm (\a_i|\a_j)}z-w)Z^{\pm}(\a_j,w)Z^{\pm}(\a_i,z),\\
&(1-w/z)_{q^{2k}}^{(\a_i|\a_j)/k}Z^+(\a_i,z)Z^-(\a_j,w)-(1-z/w)_{q^{2k}}^{(\a_i|\a_j)/k}
Z^-(\a_j,w)Z^+(\a_i,z)\\
&=\frac{\delta_{ij}}{q_i-q_i^{-1}}\{k_i\delta(q^{-k}z/w)-
k^{-1}_i\delta(q^kz/w)\}
\endalign
$$
\endproclaim
\demo{Proof} The relations are consequences of the relations satisfied
by the operators $E^{\pm}_-(z)$ and $E^{\pm}_+(z)$. For details, see
\cite{J2}
\hfill$\square$
\enddemo 

We will call the algebra generated by operators satisfying the relations of
Theorem 2.2 
as the quantum $Z$-algebra. 
We may also require that the relations contain the
 Serre relations (which are derived from the Serre 
relations of $X_i^{\pm}(z)$) to ensure that it will come from a representation
of the quantum affine algebra. 
We define the vacuum space $\Omega_V$ as an invariant subspace 
of the quantum $Z$-operators.
$$
\Omega_V=\{ v\in V | E_+^{\pm}(\a_i, z).v=E_+^{\pm}(-\a_i, z).v=v\}.
$$

As in the simply laced case \cite{J2} we have 

\proclaim{Theorem 2.3} Let $V$ be a level $k$ module, then we have 
a linear
isomorphism:
$$\align
K(k)\otimes_{\Bbb C} \Omega_V \longrightarrow  & V\\
u\otimes w \longmapsto & i(u)\cdot w
\endalign
$$
where we identify the space $K(k)$ with the Heisenberg algebra
via the classical action, and the inclusion $i$ embeds 
$U_q(\widehat{\goth h}^-)$ into $\alg$. Subsequently we can write
$$x^{\pm}_i(z)=E^{\pm}_-(-\a_i, z)E^{\pm}_+(-\a_i, z)\otimes Z^{\pm}(\a_i, z).
$$
\endproclaim
\demo{Proof} The proof is exactly the same as in the simply laced case,
as we observe that $Z^{\pm}(\a_i, z)$ commutes with the operators
$a_i(n)$ of the Heisenberg subalgebra.
\enddemo\hfill$\square$

Thus the construction of irreducible modules is reduced to that of 
corresponding $Z$-vertex operators. 

\proclaim{Proposition 2.4}
A $\alg$-module $V$ is irreducible if and only if
the vacuum space $\Omega_V$ is an irreducible $U_q(\widehat{h})$-module.
\endproclaim
\hfill$\square$

\head 3. Various constructions of $\alg$
\endhead

\subheading{3.1} $k=1$, Frenkel-Jing realization for $(ADE)^{(1)}$ \cite{FJ}.

Let $\Bbb C\{Q\}$ be the twisted group algebra of $Q$ 
generated by $e^{\a}$, $\a\in Q$ subject to the relations
$$
e^{\a}e^{\be}=(-1)^{(\alpha|\beta)}e^{\be}e^{\a}.
$$
For each level one integrable module $\Lambda_r$, there corresponds
to a minuscule weight $\lambda_r\in Q$ (cf. Table 1.1). We formally adjoin $e^{\la_r}$
to form the vector space $\Bbb C\{Q\}e^{\lambda_r}$.
 
We define the action of $z^{\partial_{\a}}$ on  $\Bbb C\{Q\}e^{\lambda_r}$
by 
$$
z^{\partial_{\a}}e^{\be}e^{\la_r}=z^{(\a|\be+\omega_r)}e^{\be}e^{\la_r}
$$ 
for $\a\in Q$.

The following construction is a reformulation of the 
basic vertex representations given in \cite{J2}. 

The level one module $V(\Lambda_r)$ is realized as
$$
V=Sym(\widehat{\goth h}^-)\otimes \Bbb C\{Q\}e^{\lambda_r}
$$
with the vacuum space equals to $\Bbb C\{Q\}e^{\lambda_r}$ and the $Z$-operators
are given as:
$$
Z^{\pm}(\a_i)=e^{\pm\a_i}z^{\pm \partial_{\a_i}},
$$
and the highest weight vector is $e^{\la_r}$.

The construction was generalized to all twisted types in \cite{J1}.

\subheading{3.2} $B_l^{(1)}$, k=1, \cite{B, JM}.

There are three level one modules: $V(\Lambda_0)$, $V(\Lambda_1)$, and 
$V(\Lambda_l)$. We can write $\Lambda_i=\Lambda_0+\la_i$, where $\la_i$ are
the fundamental weights for $\goth g$.

Let $\Bbb C\{Q\}$ be the twisted group algebra generated by $e^{\a}$, $\a\in Q$
subject to the relations:
$$
e^{\a}e^{\be}=(-1)^{(\a|\be)+(\a|\a)(\be|\be)}
e^{\be}e^{\a}.
$$
The vector space $\Bbb C\{P\}$ is defined as in \S 3.1. Let $Q_0$ be the 
sublattice of the long roots, then we have
$$\Bbb C\{Q\}=\Bbb C\{Q_0\}\oplus \Bbb C\{Q_0\}e^{\la_1}.
$$

We set $Z=\Bbb Z+s$, $s=0$ or $1/2$. The $q$-deformed Clifford algebra
$C_q^s$ is the quadratic algebra generated by $\kappa(k)$ satisfying the
relations:
$$
\{\kappa(m), \kappa(n)\}=(q^{m}+q^{-m})\delta_{m, -n},  \tag3.1
$$
where $\{a, b\}$ is the anticommutator, and $m, n\in Z$.

Let $\Lambda(C^-_q)^s$ be the polynomial algebra generated by $\ka(-k)$, 
$k\in Z_{>0}$, and $\Lambda(C^-_q)^s_0$ (resp. $\Lambda(C^-_q)^s_1$) be the 
subalgebra 
generated by products of even (resp. odd) number of generators $\ka(-k)$'s. 
Then
$$
\Lambda(C^-_q)^s=\Lambda(C^-_q)^s_0\oplus \Lambda(C^-_q)^s_1.
$$

The algebra $C_q^{s}$ acts on $\Lambda(C^-_q)^s$ canonically by 
($k\in Z_{>0}$)
$$\align
\kappa(-n).v &= \kappa(-n)\wedge v,  \qquad\qquad\text{for $n>0$}\\
\kappa(n).1  &= 0 \qquad\qquad\qquad\text{for $n>0$}
\endalign
$$ 

We define that
$$\align
W_0&=Sym(\hat{\goth h}^-)\otimes ({\Bbb C}\{Q_0\}
\otimes\Lambda(C^-_q)^{1/2}_0\oplus {\Bbb C}\{Q_0\}e^{\lambda_1}
\otimes\Lambda(C^-_q)^{1/2}_1), \\
W_1&=Sym(\hat{\goth h}^-)\otimes ({\Bbb C}\{Q_0\}e^{\lambda_1}
\otimes\Lambda(C^-_q)^{1/2}_0\oplus {\Bbb C}\{Q_0\}
\otimes\Lambda(C^-_q)^{1/2}_1),\\
W_l&=Sym(\hat{\goth h}^-)\otimes {\Bbb C}\{Q\}e^{\la_l}\otimes \Lambda(C^-_q)^0.
\endalign
$$
Then the algebras $U_q(\hat{\goth h})$, ${\Bbb C}\{P\}$ and $C_q^s$ (for suitable
$s=0, 1/2$) act on the spaces $W_i$ canonically by extending their respective
actions component-wise. We also define the degree action by
$$
d.f\otimes e^{\a}\otimes g=(-\sum_{j=1}^km_j-\sum_{j=1}^ln_j-
\frac{(\be|\be)}2 +\frac{(\alpha_i|\alpha_i)}2 )f\otimes e^{\a}\otimes g,
$$
where
$f\otimes e^{\a}\otimes g=
a_{i_1}(-m_{1})\cdots a_{i_k}(-m_k)\otimes e^{\a}\otimes
\ka(-n_1)\cdots \ka(-n_l)
\in W_i$. For convenience we denote $\la_0=0$.

\proclaim{Proposition 3.2}  
The space $W_i$ ($i=0, 1, l$) is  the irreducible representation 
$V(\Lambda_i)$  of the quantum affine Lie algebra
$U_q(B_l^{(1)})$ under the action:
$$
\align
 \gamma &\mapsto  q, \ 
K_j \mapsto q^{\partial_{\a_j}},  \  a_j(k)\mapsto a_j(k) 
\quad (1\leq j\leq l),\\
Z_i^{\pm}(z)&=e^{\pm\alpha_i}z^{\pm\partial_i} , \qquad i=1, \cdots, l-1\\
Z_l^{\pm}(z)&=e^{\pm\alpha_l}z^{\pm\partial_l}(\pm\ka(z)),
\endalign
$$
where $\ka(z)=\sum_{m\in Z}\ka(m)z^{-m}$ 
and the degree operator $d$ is defined above. 
The highest weight vectors are respectively: 
$$
 |\Lambda_i\rangle
=1\otimes e^{\lambda_i}\otimes 1, \quad i=0, 1, l. 
$$
\endproclaim

\subheading {3.3} $k=2$, $A_1^{(1)}$ \cite{I, J2}

There are three level two integrable modules $V(2\Lambda_0)$, $V(2\Lambda_1)$,
and $V(2\Lambda_0+2\Lambda_1)$. 

We need the Clifford algebra $C_{q^2}^{1/2}$ (3.1) to realize level $2$ 
$U_q(\hat{sl}_2)$. 

Define the action of $d$ on $\Gamma(\Cal C_{q^2}^s)$ by
$$
d.\kappa(n_1)\wedge \cdots \wedge\kappa(n_m)=(n_1+\cdots +n_m)
\kappa(n_1)\wedge \cdots \wedge\kappa(n_m).
$$

\proclaim{Proposition 3.3} \cite{J2} Set 
$$
\Omega=\Gamma(\Cal C_{q^2}^{-, 1/2})\otimes \Bbb C[Q],
$$
and extend naturally the actions of $e^{\beta}$ and $d$ to the space $\Omega$.
We can define a representation of $U_q(\widehat{sl}_2)$ on $V=K(2)\otimes \Omega$
via:
$$
Z^{\pm}(z)=\pm \kappa(z)\otimes e^{\pm\a}z^{\pm \a/2+1/2}.
$$ 
Moreover, the modules $V$ decomposes itself as
$$
V=L(2\Lambda_0)\oplus L(2\Lambda_1-(1/2)\delta).
$$
\endproclaim

For the level two module $L(\Lambda_0+\Lambda_1)$ we consider another Clifford
algebra $\Cal C_{q^2}^0$. 

In order to construct the vacuum space for the representation
$L(\Lambda_0+\Lambda_1)$ we need to extend the space 
$\Gamma(\Cal C_{q^2}^{-,0})$.
Define the action of $\kappa(0)$ on the two-dimensional space $\Bbb C^2$ by
$$
\kappa(0).\pmatrix u\\v\endpmatrix=\pmatrix v\\u\endpmatrix
$$
and then naturally extend the action of $\Cal C$ on the space 
$\Gamma(\Cal C_{q^2}^{-,0})\otimes \Bbb C^2$.

\proclaim{Proposition 3.4} \cite{J2}
Let $V=K(2)\otimes \Omega_V$ with
$$
\Omega_V=\Gamma(\Cal C_{q^2}^{-,0})\otimes \Bbb C^2\otimes \Bbb C[Q]e^{\a/2},
$$ 
then we can realize a representation of $U_q(\widehat{sl}_2)$ by
$$
Z^{\pm}(z)=\pm \kappa(z)\otimes e^{\pm \a}z^{\pm\a/2+1/2}.
$$
Moreover, we have
$$
V=L(\Lambda_0+\Lambda_1-\delta/8)\oplus L(\Lambda_0+\Lambda_1-\delta/8).
$$
\endproclaim
\demo{Proof} The proposition is proved similarly as Proposition 3.2.
\hfill$\square$
\enddemo

\subheading{3.4} Higher level $U_q(\widehat{sl}_2)$ 
\cite{J2, ABH, AOS, KMPY, M, S}

We consider irreducible modules of highest weights
$$
\lambda=i\Lambda_0+(k-i)\Lambda_1,
$$
where $i\in\{0, \cdots, k\}$. 
We construct the module as an submodule of the tensor product of level one 
modules $$
L(\Lambda_0)^{\otimes i}\otimes L(\Lambda_1)^{\otimes(k-i)}\tag3.2
$$
where the tensor product is induced from the Drinfeld comultiplication 
(cf. \cite{DFr} \cite{J2}).
Drinfeld's coproduct plays an important role in studying
 higher level representations based on Frenkel-Jing constructions
as we noticed in \cite{J2}. This technique was further
developed by the work \cite{DW, DFg}. 
 
In order to state Drinfeld comultiplication we introduce the following
notation. For a field operator $A(z)=\sum_n A_nz^{-n}$ acting on a space $V$
and a commutative operator $q^c$ on $V$
we define
$$\align
A(z\otimes q^c)&=\sum_n (A_n\otimes q^{nc})z^{-n}=
\sum_n A_nz^{-n}\otimes q^{nc}\\
A(q^c\otimes z)&=\sum_n (q^{nc}\otimes A_n)z^{-n}=
\sum_n q^{nc}\otimes A_n z^{-n}
\endalign
$$
as operators on $V\otimes V$. In particular, $A(z\otimes 1)=A(z)\otimes 1$.
The Drinfeld comultiplication is given by the following formula:
$$\align
\Delta(\phi_i(z))&=\phi_i(z\otimes q^{-c/2}) 
\phi_i(q^{c/2}\otimes z), \tag3.3 \\
\Delta(\psi_i(z))&=\psi_i(z\otimes q^{c/2}) 
\psi_i(q^{-c/2}\otimes z), \tag3.4\\
\Delta(x^+_i(z)) &=x^+_i(z\otimes 1)+\psi_i(q^{-c/2}z\otimes 1)
           x^+_i(q^{-c}\otimes z), \tag3.5\\
\Delta(x^-_i(z)) &=x^-_i(1\otimes z)+x^-_i(z\otimes q^{-c}) 
\phi_i(1\otimes q^{-c/2}z). \tag3.6
\endalign
$$
Thus the  first formula means that 
$$\align
\Delta(\phi_i(z))&=\sum_{m,n}(\phi_i(-n)z^n\otimes q^{-nc/2})\cdot
                            (q^{mc/2}\otimes \phi_i(-m)z^m)\\
&=\sum_{n=0}^{\infty}\sum_{m=0}^{n}\phi_i(-m)q^{(n-m)c/2}\otimes
           \phi_i(-n+m)q^{-mc/2}z^n.
\endalign
$$

The counit $\epsilon$ and antipode $S$ are given by
$$\align
\epsilon(1)&=\epsilon(\phi_i(z))=\epsilon(\psi_i(z))=1,\\
\epsilon(x_i^{\pm}(z))&=0, \\
S(\phi_i(z))&=\phi_i(z)^{-1}, \qquad S(\psi_i(z))=\psi_i(z)^{-1}, \\
S(x_i^+(z))&=-\psi_i(q^{c/2}z)^{-1}x_i^+(q^cz)\\
S(x_i^-(z))&=-x_i^-(q^cz)\phi_i(q^{c/2}z)^{-1}.
\endalign
$$

In \cite{J2}
we showed that the comultiplication is equivalent to the following
component form:
$$\align
\Delta(\a_i(n)) &=q^{|n|c/2}\otimes \a_i(n) +\a_i(n)\otimes q^{-|n|c/2}, 
\tag3.7\\
\Delta(x^+_i(n)) &=x^+_i(n)\otimes 1+\sum_{m\geq 0}\psi_i(m)q^{(m+2n)c/2}\otimes
                                 x^+_i(n-m), \tag3.8\\
\Delta(x^-_i(n)) &=1\otimes x^-_i(n)+\sum_{m\geq 0}x^-_i(n+m)\otimes\phi_i(-m)
q^{(m+2n)c/2}.\tag3.9
\endalign
$$

We extend the comultiplication to multi-tensor product in the natural way.
Define
$$
\Delta_{k-1}: \U \longrightarrow \U^{\otimes k} 
$$
inductively
$$
\Delta_1=\Delta, \Delta_j=(1\otimes \Delta)\Delta_{j-1}=(\Delta\otimes 1)\Delta_{j-1}, 
\qquad\text{for}\quad j=2, \cdots, k-1
$$

For $\e_1, \cdots, \e_s\in\{\pm\}=\{\pm1\}$ we have
$$\align
X^{\e_1}_1(z_1)\cdots X^{\e_s}_s(z_s)&=E^{\e_1}_-(-\a_1, z_1)\cdots
E^{\e_s}_-(-\a_s, z_s)\\
&E^{\e_1}_+(-\a_1, z_1)\cdots E^{\e_s}_+(-\a_s, z_s)
\Z
\endalign
$$
where the-multi argument $Z$-operator $\Z$ is defined by
$$
\Z=\prod_{1\leq i<j\leq s}(1-\frac{z_j}{z_i}
q^{-(\e_i+\e_j)k/2})^{-2\e_i\e_j/k}_{q^{2k}}
Z^{\e_1}(\a_1, z_1)\cdots Z^{\e_s}(\a_s, z_s).
$$

In particular, we have 
$$
Z(\epsilon, z)=Z^{\epsilon}(\a, z).
$$

We define the operators $Z(\epsilon_1\cdots\epsilon_s;n_1\cdots n_s)$ as the  
components of the Z-opera\-tors: 
$$
\Z=\sum Z(\epsilon_1\cdots\epsilon_s;n_1\cdots n_s)z^{-n_1}\cdots z^{-n_s}.
$$
Clearly the operators $\W$ preserve $\Omega_V$, and 
$$
deg\, \W=n_1+\cdots+n_s .
$$

Using the results of poles \cite{DW} and Drinfeld comultiplication, Ding and
Feigin \cite{DFg} recently gave a comprehensive study of higher level
representations of the quantum affine algebra $U_q(\widehat{sl_2})$.

In the above we consider examples of irreducible modules. We now turn to
reducible modules. According to Proposition 2.4 we only need to construct
quantum $Z$-algebras acting on $U_q(\widehat{\goth h})$-reducible 
modules. Usually we need to enlarge
the algebra by adding auxiliary operators.
The quantum Wakimoto constructions are done 
in this way. In the following we will only discuss one
of the constructions, since all are closely related \cite{Bo}.
Level $k$ modules of $\widehat{sl_2}$ are also studied in \cite{KMPY}
using crystal graphs.

Consider three independent 
sets of bosonic operators $a(n), b(n), c(n)$ satisfying the
commutation relations:
$$\align
[a(n), a(m)]&=\frac{[2n]}{n}[nk]\delta_{n, -m}\\
[b(n), b(m)]&=-\frac{[2n]}{n}[nk]\delta_{n, -m}\\
[c(n), c(m)]&=\frac{[2n]}{n}[n(k+2)]\delta_{n, -m}
\endalign
$$
and we let $\beta$ be an auxiliary root such that $(\beta|\beta)=
-(\alpha|\alpha)$ and $(\alpha|\beta)=0$.

\proclaim{Proposition 3.5} The realization space of $U_q(A_1^{(1)})$ 
at level $k$ is defined as 
$$
V=Sym(a(-n)'s, b(-n)'s, c(-n)'s)\otimes C[\Bbb Z{\a}+\Bbb Z{\beta}].
$$
and the action is given by
$$\align
Z^{\pm}(z) &\longmapsto
             e^{\pm \a_1}z^{\pm\partial_{\alpha}}
             \frac {U^{\pm}(z)}
             {q-q^{-1}}(B_+(zq^{\mp (k+2)/2}z)W_+(q^{\mp k/2}z)^{\pm 1}\\
           &\qquad\qquad -B_-(zq^{\pm (k+2)/2}z)W_-(q^{\pm k/2}z)^{\pm 1}),
\endalign
$$
where 
$$\align
U^{\pm}(z)&=exp(\pm\sum_{m=1}^{\infty}\frac{q^{\mp km/2}b(-m)}{[km]}z^{m})
exp(\mp\sum_{m=1}^{\infty}\frac{q^{\mp km/2}b(m)}{[km]}z^{-m})e^{\pm \beta}
z^{\pm \partial_{\beta}}\\
B_{\pm}(z)&=exp(\pm(q-q^{-1})\sum_{m=1}^{\infty}\frac{[m]}{[2m]}b(\pm)
z^{\mp m})q^{\mp \partial_{\beta}}\\
W_{\pm}(z)&=exp(\pm(q-q^{-1})\sum_{m=1}^{\infty}\frac{[m]}{[2m]}c(\pm)
z^{\mp m}).
\endalign
$$
\endproclaim

\subheading{3.5} Realizations of $U_q(C_l^{(1)})$
and $U_q(F_4^{(1)})$.

Recently we have constructed 
the quantum symplectic affine algebra $U_q(C_n^{(1)})$ 
at level $-1/2$ \cite{JKM1} and level $1$ \cite{JKM2}.
The reducible modules were constructed using
 the same idea of the quantum 
$Z$-algebras. In the case of level one our realization contains 
all integrable level one modules. 

Recalling the construction in 3.1 we also formally define the space
$C[P]$ as the direct sum $C[Q]\oplus C[Q]e^{\lambda_1}$ by adjoining 
the element $e^{\lambda_1}$. The short roots $Q_s$ constitute the root system
$A_{n-1}$. Let $\tilde{Q}$ be the root lattice of $A_{l-1}$ generated
$\tilde{\a_i}$, $i=1, \cdots, l-1$ and normalized
so that $(\tilde{\a}_i|\tilde{\a_i})=1$, $i=1, \cdots, l-1$.  

Let $a_i(n)$ $(i=1, \cdots, l)$ and $b_i(n)$ $(i=1, \cdots, l-1)$
be two set of bosonic operators satisfying
$$\align
[a_i(n), a_j(m)]&=\frac{[(\a_i|\a_j)n]}{n}[n]\delta_{n, -m},\\
[b_i(n), b_j(m)]&=\frac{[(\a_i|\a_j)n]}{n}[n]\delta_{n, -m}.
\endalign
$$

Then level one realization is defined on the Fock space 
$$V=Sym(a_i(-n)'s, b_i(-n)'s)\otimes C[P]\otimes C[\tilde{Q}].$$

The action of the $Z$-operators are defined by
$$\align
Z_i^{\pm}(z)&=(U_i^+(zq^{-1/2})+U^-_i(zq^{1/2})(-1)^{2\partial_i})e^{\pm \a_i}
z^{\pm \partial_i}, \qquad\qquad i=1, \cdots l-1,\\
Z_l^{\pm}(z)&=e^{\pm \a_l}z^{\pm \partial_l}.
\endalign
$$
where the auxiliary operator $U^+_i(z)$ and $U_i^-(z)$ are defined by:
$$
U_i^{\pm}(z)=exp(\pm\sum_{n\geq 1}\frac{b_i(-n)}{[n]}z^n)
exp(\mp\sum_{n\geq 1}\frac{b_i(-n)}{[n]}z^{-n})e^{\pm\tilde{\a_i}}
z^{\pm \partial_i}
$$
and the operator $e^{\pm\tilde{\a}_i}$ acts on $C[A_{l-1}]$ (the group algebra
generated by the short roots) by
$$
e^{\tilde{\a_i}}.e^{\tilde{\beta}}=\epsilon(\a_i, \beta)
e^{\tilde{\a_i}+\tilde{\beta}}.
$$
where the function $\epsilon(\a_i, \a_j)$ is given by
$$
\epsilon(\a_i, \a_j)=\cases -1, & i=j\\
                             1, & i>j\\
                        (-1)^{a_{ij}}, & i<j\endcases
$$

To construct the quantum affine algebra $U_q(F_4^{(1)})$, we let 
$b_i(n)$ ($i=1, 2$) be the auxiliary bosonic operators satisfying the 
same relations as the first two Heisenberg generators associated to
$U_q(F_4^{(1)})$. Now the short roots form a sublattice
isomorphic to that of $A_2$. The representation space is defined as
$$
V=Sym(a_i(-n)'s, b_j(-n)'s)\otimes C[P]\otimes C[A_{2}].
$$
where $i=1, 2$ and $j=1, 2$. 

Similar to the case of $C_n^{(1)})$ we have the following result.
Note that we do not need to adjust signs by group algebra.

\proclaim{Theorem 3. 6} For $i=1, \cdots 4$ and $j=1, 2$, The following action defines 
a level one representation of $U_q(F_4^{(1)})$
$$\align
Z_i^{\pm}(z)&=(U_i^+(zq^{-1/2})+U^-_i(zq^{1/2})(-1)^{2\partial_i})e^{\pm \a_i}
z^{\pm \partial_i}, \qquad\qquad i=1, 2,\\
Z_j^{\pm}(z)&=e^{\pm \a_j}z^{\pm \partial_j}. \qquad\qquad i=3,4,
\endalign
$$
where the auxiliary operator $U^+_i(z)$ and $U_i^-(z)$ are defined by the 
same expressions as in the case of $C_n^{(1)}$.
\endproclaim

It is clear that the operators $U_i^{\pm}(z)e^{\mp b_i}z^{\mp b_i(0)}$ 
satisfy the relations of quantum $Z$-algebras associated with the
the same invariant form $(\ \ |\ \ )$. Thus they are some auxiliary
quantum $Z$-operators in the representation. This phenomena appear in 
almost all (vertex operator) reducible modules. To get irreducible
modules one can use quantum screening operators \cite{M}.

We conclude this subsection by a level $-1/2$ representation of 
$U_q(C_l^{(1)})$ \cite{JKM1}.

Let $a_i(m)$ and $b_i(m)$ $(i=1, \cdots, l; m, n\in \Bbb Z)$
be the operators satisfying the following relations:
$$\align
[a_i(m), a_j(n)]&=\delta_{m+n,0}\frac{[m(\alpha_i|\alpha_j)][-m/2]}m,\\
[b_i(m), b_j(n)]&=m\delta_{ij}\delta_{m+n,0}\\
[a_i(m), b_j(n)]&=0.
\endalign
$$  
 The module ${\Cal F}_{\a,\be}={\Cal F}^1_{\a}\otimes
{\Cal F}^2_{\be}$ is built from the Heisenburg algebra
   by the defining relations
   $$ a_i(m) | \a, \be \rangle = 0 \quad (m>0)
      \; , \quad
      b_i(m) | \a, \be \rangle = 0 \quad (m>0) \;, $$
   $$ a_i(0) | \a, \be \rangle = (\a_i|\a) | \a, \be \rangle
      \; , \quad
      b_i(0) | \a, \be \rangle = (2\vep_i|\be) | \a, \be \rangle \; , $$
   where $| \a, \be \rangle=| \a\rangle\otimes |\be \rangle$
 $(\a\in  P+\frac{{\bold Z}}2\la_l, \be\in P)$ is  
   the vacuum vector. 
   The grading operator ${d}$ is defined by
 $$ d. | \alpha,\beta \rangle
      =( (\alpha|\alpha)-(\beta|\beta-\lambda_n))
       |\alpha,\beta \rangle . $$
   
   We set the Fock space
   $$ \widetilde{\Cal F}:=
          \bigoplus_{\a\in P+\frac12{\bold Z}\la_l,\be \in P}
          {\Cal F}_{\a, \be} . $$
   Let $e^{a_i}$ and $e^{b_i}$ be operators
   on ${\widetilde{\Cal F}}$ given by:
   $$ e^{a_i}|\a,\be \rangle= |\a+\a_i,\be \rangle
      \quad , \quad
      e^{b_i}|\a, \be \rangle= |\a,\be+\vep_i \rangle . $$
   Let $: \quad :$ be the usual bosonic normal ordering defined by
   $$ :a_i(m) a_j(n): = a_i(m) a_j(n) \, (m \leq n) , 
   \; a_j(n) a_i(m) \, (m>n) , $$
   $$ :e^{a} a_i(0):=:a_i(0) e^{a}:=e^{a} a_i(0) \, , $$
   $$ :e^{b} b_i(0):=:b_i(0) e^{b}:=e^{b} b_i(0) \, . $$
   and similar normal products for the $b_i(m)$'s. 
   
Let $\partial=\partial_{q^{1/2}}$ be the $q$-difference operator:
$$\partial_{q^{1/2}}(f(z))=\frac{f(q^{1/2}z)-f(q^{-1/2}z)}
{(q^{1/2}-q^{-1/2})z}
$$
   We introduce the following operator:
$$
U^{\pm}_i(z)=
             \exp ( \pm \sum^{\infty}_{k=1} \frac{b_i(-k)}{k} z^k)
             \exp ( \mp \sum^{\infty}_{k=1} \frac{b_i(k)}{k} z^{-k})
             e^{\pm b_i} z^{\pm b_i(0)} .
$$
\proclaim{Proposition 3.7}The Fock space
    $\widetilde{\Cal F}$ is a $U_q$-module of level $-\frac{1}{2}$
    under the action defined by 
    $ \gamma \longmapsto q^{-\frac{1}{2}},
       K_i \longmapsto q^{a_i(0)}$, 
    $ a_{im} \longmapsto a_i(m), 
       q^d \longmapsto q^{{d}}$, and 
$$\align 
Z_i^+(z) &\longmapsto
             \partial U^+_i(z)U^-_{i+1}(z)e^{\alpha_i}z^{-2a_i(0)}, \qquad i=1, \cdots n-1\\
Z_i^-(z) &\longmapsto
             U^-_i(z)\partial U^+_{i+1}(z)e^{-\alpha_i}z^{2a_i(0)}, \qquad i=1, \cdots n-1\\
Z_l^+(z) &\longmapsto
             \left(
             \frac{1}{q^{\frac{1}{2}}+ q^{-\frac{1}{2}}}
             :U^+_l(z) \partial^2 U^+_l(z):
     -:\partial U^+_l(q^{-\frac{1}{2}}z)\partial U_l^+(q^{-\frac12}z):
             \right)
            e^{\alpha_l}z^{-2a_l(0)} \\
Z_n^-(z) &\longmapsto
             \frac{1}
                  {q^{\frac{1}{2}}+ q^{-\frac{1}{2}}}
             :U^-_l(q^{\frac{1}{2}}z)U^-_l(q^{-\frac{1}{2}}z):
         e^{-\alpha_l}z^{2a_l(0)}.
\endalign
$$
\endproclaim


\subheading{3.7} Level one realization of $U_q(G_2^{(1)})$ \cite{J3}.

We can realize level one $U_q(G_2^{(1)})$ by introducing two sets of
auxiliary bosonic operators $b_2(n)$ and $c_2(n)$ ($n\in \Bbb Z$) with the 
defining relations:
$$\align
[b_2(n), b_2(m)]&=-\frac{[\frac23n]}{n}[n]\delta_{n, -m}\\
[c_2(n), c_2(m)]&=\frac{[\frac23n]}{n}[\frac53n]\delta_{n, -m}\\
[a_i(n), b_2(n)]&=[a_i(n), c_2(n)]=0
\endalign
$$

Let us enlarge the root lattice by adding
an auxiliary root $\beta_2$ with $(\be_2|\be_2)=-2/3$. The actions of
$a_i(m), b_2(m)$ are defined similarly as previous sections. We often
identify $a_i(0)$ with $\partial_{\alpha_i}$.

\proclaim{Proposition 3.8} On the Fock space 
$$V=Sym(a_i(-n), b_2(-n), c_2(-n))\otimes C[P+\Bbb Z\beta_2]
$$
the correspondence
$$\align 
Z_1^{\pm}(z) &\longmapsto e^{\pm \a_1}z^{\pm\partial_{\alpha_1}}\\
Z_2^{\pm}(z) &\longmapsto
             e^{\pm \a_2}z^{\pm\partial_{\alpha_2}}
             \frac {U^{\pm}(z)}
             {q-q^{-1}}(B_+(zq^{\mp 5/2}z)W_+(q^{\mp 3/2}z)^{\pm 1}\\
           &\qquad\qquad -B_-(zq^{\pm 5/2}z)W_-(q^{\pm 3/2}z)^{\pm 1})
\endalign
$$
gives a level one $U_q(G_2^{(1)})$-module and where 
$$\align
U^{\pm}(z)&=exp(\pm\sum_{m=1}^{\infty}\frac{q^{\mp m/2}b_2(-m)}{[m]}z^{m})
exp(\mp\sum_{m=1}^{\infty}\frac{q^{\mp m/2}b_2(m)}{[m]}z^{-m})e^{\pm \be_2}
z^{\pm b_2(0)}\\
B_{\pm}(z)&=exp(\pm(q-q^{-1})\sum_{m=1}^{\infty}\frac{[m/3]}{[2m/3]}b_2(\pm)
z^{\mp m})q^{\mp b_2(0)/2}\\
W_{\pm}(z)&=exp(\pm(q-q^{-1})\sum_{m=1}^{\infty}\frac{[m/3]}{[2m/3]}c_2(\pm)
z^{\mp m}).
\endalign
$$
\endproclaim
The proof is obtained by a careful analysis of the Serre relations \cite{J3}.

In this paper we restricted ourselve to the untwisted types or type one.
We can similarly introduce the quantum $Z$-algebras for the twisted types
using \cite{J4}. Then all the level one modules of twisted types \cite{J1}
can be realized by the method discussed in this paper.


\widestnumber\key{KMPY}
\Refs\nofrills{References}
\ref\key ABH\by A. Abada, A.H. Bougourzi, and El. Gradechi
       \paper Deformation of Wakimoto construction
       \jour Mod. Phys. Lett A\vol 8\yr 1993\pages 715-724\endref
\ref\key AOS\by H. Awata, S. Odake, and J. Shiraishi
       \paper Free boson representation of $U_q(\hat{sl}_N)$
       \jour Commun. Math. Phys.\yr 1994\pages 61-83
       \endref
\ref\key Be\by J. Beck 
        \paper Braid group action and quantum affine algebras
        \jour Commun. Math. Phys.\vol 165\pages 555-568\yr 1994\endref
\ref\key Br\by D. Bernard 
         \paper Vertex operator representations of quantum affine algebras $U_q(B^{(1)})$
         \jour Lett. Math. Phys.\vol 17\pages 239-245
         \yr 1989\endref
\ref\key Bo\by A.H. Bougourzi
         \paper Uniqueness of the bosonization of the $U_q(su(2)_k)$ quantum current algebra
         \jour Nucl. Phys.\vol B404\yr 1993
         \pages 457-482\endref
\ref\key BV\by A.H. Bougourzi and H. Vinet
         \paper A quantum analogue of the $\Cal Z$ algebra
        \jour J. Math. Phys. \vol 37\yr 1996\pages 3548--3567
         \endref
\ref\key DFg\by J. Ding and B. Feigin
         \paper Quantum current operators. II. Difference equations of 
       quantum current operators and quantum parafermion construction
        \jour Publ. Res. Inst. Math. Sci. \vol 33\yr 1997\pages  285--300
        \endref
\ref\key DFr\by J. Ding and I.B. Frenkel
        \paper Isomorphisms of two constructions of quantum affine algebras
        \jour Comm. Math. Phys. \vol 156\yr 1993\pages 277-300\endref
\ref\key DW\by J. Ding and T. Miwa
        \paper Quantum current operators. I. Zeros and poles of quantum current
          operators and the condition of quantum integrability\jour 
         Publ. Res. Inst. Math. Sci. \vol 33\yr 1997\pages 277--284\endref
\ref\key D\by  V.G. Drinfeld
           \paper New realization of Yangians and quantum affine algebras
           \jour Sov. Math. Dokl.\vol 36\pages 212-216\endref 
\ref\key FJ\by I.B. Frenkel and N. Jing
       \paper Vertex representations of quantum affine algebras
       \jour Proc. Nat'l. Acad. USA\vol 85\pages 9373-9377\yr 1988\endref
\ref\key FLM \by I. Frenkel, J. Lepowsky, and A. Meurman\book
      Vertex operator algebras and the Monster\publ Academic Press\publaddr
       New York\yr 1988\endref 
\ref\key J1\by N. Jing
         \paper Twisted vertex representations of quantum affine algebras
         \jour Invent. Math.\vol 102\pages 663-690\yr 1990\endref
\ref\key J2\by N. Jing
         \paper Higher level representations of the quantum affine algebra $U_q(\widehat{sl}_2)$
         \jour Jour. Alg.\vol 182\pages 448-468\yr 1996\endref
\ref\key J3\by N. Jing
         \paper Level one representations of $U_q(G_2^{(1)})$\jour
         Proc. Amer. Math. Soc., to appear\endref
\ref\key J4\by N. Jing
         \paper On Drinfeld realization of quantum affine algebras
         \inbook Proc. Lie algebras and the Monster\publ W. de Gruyter
         \publaddr Berlin, to appear, q-alg/9610035\endref
\ref\key JM\by N. Jing and K.C. Misra
         \paper Vertex operators of level one $U_q(B_n^{(1)})$-modules
         \jour Lett. Math. Phys.\vol 36\yr 1996\pages 127-143\endref
\ref\key JKM1\by N. Jing, Y. Koyama, and K.C. Misra
         \paper Bosonic realizations of  $U_q(C_n^{(1)})$
         \jour Jour. Algebra \vol 200\yr 1998\pages 155-172. 
\endref
\ref\key JKM2\by N. Jing, Y. Koyama, and K.C. Misra
         \paper Level one realizations of quantum affine algebras 
         $U_q(C_n^{(1)})$
         \jour submitted (QA/9802123)\yr 1998\endref
\ref\key I\by M. Idzumi
         \paper Level two irreducible representations of $U_q(\hat{sl}_2)$
         \jour Int. J. Mod. Phys. A\vol 9\yr 1994\pages 4449-4484
         \endref
\ref\key K\by V.G. Kac
         \book Infinite dimensional Lie algebras\bookinfo 3rd ed.
         \publ Cambridge University Press\yr 1990\endref 
\ref\key KMPY\by M. Kashiwara, T. Miwa, J\.-U\.H. Peterson, and 
         C\.M Yung         
         \paper Perfect crystals and $q$-deformed Fock spaces
         \jour Selecta Mathematica\vol 2\yr 1996\pages 415--499\endref
\ref\key LP1\by J. Lepowsky and M. Primc
         \paper Standard modules for type one affine Lie algebras
         \inbook Number Theory\bookinfo  Lect. Notes in Math.
         \vol  1052\publ Springer-Verlag\publaddr New York
         \yr 1984\pages 194-251\endref
\ref\key LP2\by J. Lepowsky and M. Primc
         \book Structure of the standard modules for the affine Lie algebra $A^{(1)}_1$
         \bookinfo Contemporary Math.\vol 46\publ Amer. Math. Soc.
         \publaddr Providence, RI \yr 1985\endref
\ref\key LW1\by J. Lepowsky and R. L. Wilson
          \paper A new family of algebras underlying the Rogers-Ramanu\-jan identities and generalizations
           \jour Proc. Nat'l. Acad. Sci. USA\vol 78\pages 7254-7258\yr
           1981\endref
\ref\key LW2\by J. Lepowsky and R. L. Wilson
          \paper The structure of standard modules, I: Universal algebras and the Rogers-Ramanu\-jan identities
          \jour Invent. Math.\vol 77
          \pages 199-290\yr 1984\endref
\ref\key LW3\by J. Lepowsky and R. L. Wilson
          \paper The structure of standard modules, II: The case $A_1^{(1)}$, principal gradation
          \jour Invent. Math. \vol 77\pages 417-442
          \yr 1985\endref
\ref\key M\by A. Matsuo
       \paper A $q$-deformation of Wakimoto modules, primary fields and
       screening operators\jour Comm. Math. Phys. \vol 160 \yr 1994
       \pages 33--48\endref
\ref\key P\by M. Primc
       \paper Standard representations of $A_n^{(1)}$
       \inbook Infinite-dimensional Lie algebras and groups
       \publ World Sci. Singapore\publaddr Singapore\yr 1989\pages
       273-284\endref
\ref\key S\by J. Shiraishi
       \paper Free boson representation of $U_q(\hat{sl}_2)$
       \jour Phys. Lett. A\vol 171\yr 1992\pages 243-248
       \endref

\endRefs
\enddocument